\documentclass[A4j,11pt]{article}
\usepackage{graphics}
\usepackage{amsmath}
\usepackage{amssymb}
\usepackage{latexsym}
\usepackage{amsfonts}
\begin{document}

\title{{\bf On proper complex equifocal submanifolds}}
\author{{\bf Naoyuki Koike}}
%
%
\date{}
%
\maketitle
\begin{abstract}
First we show that a proper complex equifocal submanifold occurs as 
a principal orbit of a Hermann type action under certain condition.  
Next we show that a proper complex equifocal submanifold is an isoparametric 
submanifold with flat section under certain condition.  
\end{abstract}

\vspace{0.5truecm}

\section{Introduction}
C.L. Terng and G. Thorbergsson [TT] introduced the notion of an equifocal 
submanifold in a (Riemannian) symmetric space $G/K$, which is defined as 
a compact submanifold with trivial normal holonomy group, flat section and 
parallel focal structure.  Here the parallelity of the focal structure means 
that, for each parallel normal vector field $v$, the focal radii of 
the submanifold along the normal geodesic $\gamma_{v_x}$ 
(with $\gamma_{v_x}'(0)=v_x$) are independent of the choice of a point $x$ of 
the submanifold (with considering the multiplicities), where 
$\gamma_{v_x}'(0)$ is the velocity vector of $\gamma_{v_x}$ to $0$.  
Note that 
the focal radii of the submanifold along the normal geodesic $\gamma_{v_x}$ 
coincide with the zero points of a real valued function 
$$F_{v_x}(s):={\rm det}\left(\cos(s\sqrt{-1}\sqrt{R(v_x)})
-\frac{\sin(s\sqrt{-1}\sqrt{R(v_x)})}{\sqrt{-1}\sqrt{R(v_x)}}\circ A_{v_x}
\right)$$
over ${\Bbb R}$ defined in terms of the shape operator $A_{v_x}$ and 
the normal Jacobi operator $R(v_x)(:=R(\cdot,v_x)v_x)$, where $R$ is the 
curvaure tensor of the ambient symmetric space.  
In particular, in the case where $G/K$ is a Euclidean space, we have 
$F_{v_x}(s)={\rm det}({\rm id}-sA_{v_x})$ and hence the focal radii along 
$\gamma_{v_x}$ coincide with the inverse numbers of the eigenvalues of 
$A_{v_x}$ (i.e., the principal curvatures of direction $v$).  
Heintze-Palais-Terng-Thorbergsson [HPTT] defined a hyperpolar action on 
a symmetric space $G/K$ as a compact group action with flat section.  
Also, if $G/K$ is of compact type and if $H$ is a symmetric subgroup of $G$ 
(i.e., $({\rm Fix}\,\sigma)_0\subset H\subset{\rm Fix}\,\sigma$ for some an 
involution $\sigma$ of $G$), then they called the $H$-action on $G/K$ 
a Hermann action, where we note that Hermann actions are hyperpolar.  
They showed that principal orbits of hyperpolar actions are equifocal.  
According to the classification of hyperpolar actions on irreducible symmetric 
spaces of compact type by A. Kollross [Kol], hyperpolar actions of 
cohomogeneity greater than one on the symmetric spaces are orbit equivalent to 
Hermann actions.  O. Goertsches and G. Thorbergsson [GT] showed that principal 
orbits of a Hermann action are curvature-adapted.  Here 
the curvature-adaptedness means 
that, for each normal vector $v$ of the submanifold, the normal Jacobi 
operator $R(v)$ preserves the tangent space invariantly and that 
it commutes the shape operator $A_v$.  
U. Christ [Ch] showed that all irreducible equifocal submanifolds of 
codimension greater than one in a symmetric space of compact type are 
homogeneous.  

From these facts, we obtain the following fact.  

\vspace{0.3truecm}

\noindent
{\bf Fact 1.} {\sl All equifocal submanifolds of codimension greater than one 
in an irreducible symmetric space of compact type are principal orbits of 
Hermann actions and they are curvature-adapted.}

\vspace{0.3truecm}

When a non-compact submanifold $M$ in a symmetric space $G/K$ of non-compact 
type variates as its principal curvatures approch to zero, 
its focal set vanishes beyond the ideal boundary $(G/K)(\infty)$ of $G/K$.  
From this fact, we recognize that, for a non-compact submanifold 
in a symmetric space of non-compact type, the parallelity of the focal 
structure is not an essential condition.  
So, we ([Koi3]) introduced the notion of a complex focal radius of the 
submanifold along the normal geodesic $\gamma_{v_x}$ as the zero points of 
a complex valued function $F^{\bf c}_{v_x}$ over ${\Bbb C}$ defined by 
$$F^{\bf c}_{v_x}(z):={\rm det}\left(\cos(z\sqrt{-1}\sqrt{R(v_x)^{\bf c}})
-\frac{\sin(z\sqrt{-1}\sqrt{R(v_x)^{\bf c}})}
{\sqrt{-1}\sqrt{R(v_x)^{\bf c}}}\circ A^{\bf c}_{v_x}\right)$$
over ${\Bbb C}$, where $A^{\bf c}_{v_x}$ and $R(v_x)^{\bf c}$ are 
the complexifications of $A_{v_x}$ and $R(v_x)$, respectively.  
In the case where $M$ is of class $C^{\omega}$ (i.e., real analytic), we 
([Koi3]) defined the complexification $M^{\bf c}$ of $M$ as an 
anti-Kaehlerian submanifold in the anti-Kaehlerian symmetric space 
$G^{\bf c}/K^{\bf c}$.  
We ([Koi3]) showed that $z$ is a complex focal radius of $M$ along 
$\gamma_{v_x}$ if and only if 
$\exp^{\perp}(({\rm Re}\,z)v_x+({\rm Im}\,z)Jv_x)$ is a focal point of 
$M^{\bf c}$, where $\exp^{\perp}$ is the normal exponential map of 
$M^{\bf c}$ and $J$ is the complex structure of $G^{\bf c}/K^{\bf c}$.  
When $M$ variates as above and real analytically, 
its focal set vanishes beyond $(G/K)(\infty)$ but 
the focal set of $M^{\bf c}$ (i.e., the complex focal set of $M$) does not 
vanish.  From this fact, for non-compact submanifolds in a symmetric space 
of non-compact type, we recognize that the parallelity of the complex focal 
structure is an essential condition (even if $M$ is not of $C^{\omega}$).  
So, we [Koi2] defined the notion of 
a complex equifocal submanifold (which should be called a equi-complex focal 
submanifold precisely) as a (properly embedded) complete submanifold with 
trivial normal holonomy group, flat section and parallel complex focal 
structure.  Note that equifocal submanifolds in the symmetric space are 
complex equifocal.  In fact, since they are compact, their principal 
curvatures are not close to zero and hence the parallelity of their focal 
structure leads to that of their complex focal structure.  

In 1989, Terng [T2] introduced the notion of an isoparametric 
submanifold in a (separable) Hilbert space as a proper Fredholom submanifold 
with trivial normal holonomy group and constant principal curvatures.  
Here we note that the shape operators of the submanifold 
are compact operators and that they are simultaneously diagonalizable with 
respect to an orthonormal base.  
Also, she [T2] defined the notion of the parallel transport map 
for a compact semi-simple Lie group $G$ as a Riemannian submersion of 
a (separable) Hilbert space $H^0([0,1],\mathfrak g)$ onto $G$, where 
$H^0([0,1],\mathfrak g)$ is the space of all $L^2$-integrable paths in the Lie 
algebra $\mathfrak g$ of $G$.  
Let $G/K$ be a symmetric space of compact type, $\pi$ the natural projection 
of $G$ onto $G/K$ and $\phi$ the parallel transport map for $G$.  
Also, let $M$ be a compact submanifold in $G/K$ and $\widetilde M$ a component 
of the lifted submanifold $(\pi\circ\phi)^{-1}(M)$.  
In 1995, Terng-Thorbergsson [TT] showed that $M$ is equifocal if and only if 
$\widetilde M$ is isoparametric.  
Thus the research of an equifocal submanifold in a symmetric space of compact 
type is reduced to that of an isoparametric submanifold in a (separable) 
Hilbert space.  
By using this reducement of the research, they [TT] 
proved some facts for an equifocal submanifold in the symmetric space.  

Terng-Thorbergsson [TT] proposed the following problem:

\vspace{0.2truecm}

\noindent
{\bf Problem([TT]).} 
{\sl For equifocal submanifolds in symmetric spaces of non-compact type, 
is there a method of research similar to the above method of research by using 
the lift to a Hilbert space (for equifocal submanifolds in symmetric spaces 
of compact type)?}

\vspace{0.2truecm}

From 2002, we began to tackle to this problem.  
In 2004-2005, we [Koi2,3] constructed the similar method of research for 
complex equifocal submanifolds in symmetric spaces of non-compact type in more 
general.  We shall explain this method of research.  We first defined the 
notion of a complex isoparametric submanifold in a pseudo-Hilbert space as 
a (properly embedded) complete Fredholm submanifold with trivial normal 
holonomy group and constant complex principal curvatures.  
See the next section about the definition of a Fredholm submanifold in 
a pseudo-Hilbert space.  
Note that, for each normal vector $v$ of a complex isoparametric submanifold 
$M$, the shape operator $A_v$ of $M$ is not necessarily diagonalizable with 
respect to an orthonormal base and, furthermore, 
the complexified shape operator $A_v^{\bf c}$ also is not necessarily 
diagonalizable with respect to a pseudo-orthonormal base.  
In particular, if, for each normal vector $v$ of $M$, $A_v^{\bf c}$ is 
diagonalizable with respect to a pseudo-orthonormal base, then 
we called $M$ a {\it proper complex isoparametric submanifold}.  
Here we note that, in pseudo-Riemannian submanifold theory, 
if the complexified shape operators of a pseudo-Riemannian submanifold are 
diagonalizable with respect to a pseudo-orthonormal base, then 
the submanifold is called to be {\it proper} (see [Koi1] for example).  
Also, we [Koi2] defined the notion of the parallel transport map for 
a (not necessarily compact) semi-simple Lie group $G$ as a pseudo-Riemannian 
submersion of a pseudo-Hilbert space onto $G$.  
Let $G/K$ be a symmetric space of non-compact type, $\pi$ the natural 
projection of $G$ onto $G/K$ and $\phi$ the parallel transport map for $G$.  
Also, let $M$ be a (properly embedded) complete submanifold in $G/K$ and 
$\widetilde M$ a component of the lifted submanifold 
$(\pi\circ\phi)^{-1}(M)$.  
We [Koi2] showed that $M$ is complex equifocal if and only if 
$\widetilde M$ is complex isoparametric.  
Thus the research of complex equifocal submanifolds in symmetric spaces of 
non-compact type is reduced to that of complex isoparametric submanifolds in 
pseudo-Hilbert spaces.  
If each component of $(\pi\circ\phi)^{-1}(M)$ is proper complex isoparametric, 
then we ([Koi4]) called $M$ a {\it proper complex equifocal submanifold}.  
Since the shape operators of a proper complex isoparametric submanifold is 
simultaneously diagonalizeble with respect to a pseudo-orthonormal base, the 
complex focal set of the submanifold at any point $x$ consists of 
infinitely many complex hyperplanes in the complexified normal space at $x$ 
and the group generated by the complex reflection of order two with respect to 
the complex hyperplanes is discrete.  
From this fact, it follows that the same fact holds for the complex focal set 
of a proper complex equifocal submanifold.  

Let $G/K$ be a symmetric space of non-compact type and $H$ 
a symmetric subgroup of $G$ (i.e., 
$({\rm Fix}\,\sigma)_0\subset H\subset{\rm Fix}\,\sigma$ for some involution 
$\sigma$ of $G$).  Then the $H$-action on $G/K$ is called a 
{\it Hermann type action}.  
We [Koi4] showed that principal orbits of a Hermann type action on a symmetric 
space of non-compact type are proper complex equifocal and curvature-adapted.  

In 2006, Heintze-Liu-Olmos [HLO] defined the notion of isoparametric 
submanifold with flat section in a general Riemannian manifold as 
a (properly embedded) complete submanifold with flat section and trivial 
normal holonomy group whose sufficiently close parallel submanifolds have 
constant mean curvature with respect to the radial direction.  
For a compact submanifold with trivial holonomy group and flat section in a 
symmetric space of compact type, they [HLO] showed that it is equifocal 
if and only if, for each parallel normal vector field $v$, $F_{v_x}$ is 
independent of the choice of a point $x$ of the submanifold, where $F_{v_x}$ 
is the function defined in Page 1.  Thus if it is an isoparametric submanifold 
with flat section, then it is equifocal.  
Furthermore, for a compact submanifold in a symmetric space of compact type, 
they [HLO] showed that it is equifocal if and only if it is an isoparametric 
submanifold with flat section.  
On the other hand, we [Koi3] showed that, 
for a (properly embedded) complete submanifold with trivial normal holonomy 
group and flat section in a symmetric space of non-compact type, it is an 
isoparametric submanifold with flat section if and only if, for each parallel 
normal vector field $v$, $F^{\bf c}_{v_x}$ is independent of the choice of 
a point $x$ of the submanifold, where $F^{\bf c}_{v_x}$ 
is the function defined in Page 2.  Thus if it is an isoparametric submanifold 
with flat section, then it is complex equifocal.  
Furthermore, we [Koi3] showed that, if it is curvature-adapted and complex 
equifocal, then it is an isoparametric submanifold with flat section.  

For a submanifold $M$ in a Hadamard manifold $N$, we ([Koi11]) defined the 
notion of a focal point of non-Euclidean type on the ideal boundary 
$N(\infty)$.  
We [Koi11] showed that, for a curvature-adapted complex equifocal 
submanifold $M$, it is proper complex equifocal if and only if it admits 
no non-Euclidean type focal point on the ideal boundary of the ambient 
symmetric space.  
According to Theorems A and C (also Remark 1.1) in [Koi6], 
it is shown that any irreducible homogeneous complex equifocal submanifold 
of codimension greater than one admitting a totally geodesic focal 
submanifold (or a totally geodesic parallel submanifold) occurs as 
a principal orbit of a Hermann type action.  
On the other hand, we [Koi10] showed the following homogeneity theorem:

\vspace{0.3truecm}

{\sl All irreducible proper complex equifocal $C^{\omega}$-submanifolds of 
codimension greater 

than one are homogeneous.}

\vspace{0.3truecm}

\noindent
{\bf Assumption.} In the sequel, we assume that all submanifolds are 
of class $C^{\omega}$.  

\vspace{0.3truecm}

In this paper, we first prove the following fact in terms of these facts.  

\vspace{0.5truecm}

\noindent
{\bf Theorem A.} {\sl All irreducible curvature-adapted proper complex 
equifocal submanifolds of codimension greater than one in a symmetric 
space of non-compact type occur as principal orbits of Hermann type actions on 
the symmetric space.}

\vspace{0.5truecm}

\noindent
{\it Remark 1.1.} In this theorem, we cannot replace "proper complex 
equifocal" to "complex equifocal".  In fact, principal orbits of the 
$N$-action on an irreducible symmetric space $G/K$ of non-compact type and 
rank greater than one are irreducible curvature-adapted complex equifocal 
submanifolds of codimension greater than one but they do not occur as 
principal orbits of a Hermann type action, where $N$ is the nilpotent part in 
the Iwasawa's decomposition $G=KAN$ of $G$.  

\vspace{0.3truecm}

Next we prove the following fact.  

\newpage


\noindent
{\bf Theorem B.} {\sl Let $G/K$ be a symmetric space of non-compact type and 
rank $r$.  Then all proper complex equifocal submanifolds of codimension $r$ 
in $G/K$ are isoparametric submanifolds with flat section.}

\vspace{0.3truecm}

Also, we prove the following fact.  

\vspace{0.3truecm}

\noindent
{\bf Theorem C.} {\sl Let $G/K$ be a symmetric space of non-compact type and 
rank $r$ whose root system is reduced.  
Then all proper complex equifocal submanifolds of codimension $r$ in $G/K$ 
are curvature-adapted.  
}

\vspace{0.3truecm}

\noindent
{\it Remark 1.2.} By imitating the discussion in the proof of Theorem C, 
we can show the following fact:

\vspace{0.15truecm}

{\sl  Let $G/K$ be a symmetric space of compact type and rank $r$ whose 
root system is reduced.  Then all equifocal submanifolds of codimension $r$ 
in $G/K$ are curvature-adapted.}

\vspace{0.3truecm}

From Theorems A and C, the following fact follows directly.  

\vspace{0.3truecm}

\noindent
{\bf Theorem D.} {\sl Let $G/K$ be a symmetric space of non-compact type and 
rank $r(\geq2)$ whose root system is reduced.  
Then all irreducible proper complex equifocal submanifolds of codimension $r$ 
in $G/K$ occurs as principal orbits of Hermann type actions on $G/K$.}

\section{Basic notions and facts}
In this section, we recall the notions of a complex equifocal submanifold, 
the parallel transport map for a semi-simple Lie group and 
a proper complex isoparametric submanifold in a pseudo-Hilbert space.  
Since these notions are not well-known for the experts of this topic, 
we explain them in detail.  
We first recall the notion of a complex equifocal submanifold.  
Let $M$ be a (properly embedded) complete submanifold in a symmetric space 
$N=G/K$ of non-compact type.  
Assume that $M$ has flat section, that is, $\exp^{\perp}(T_x^{\perp}M)$ is 
a flat totally geodesic submanifold in $N$ for each $x\in M$, where 
$T^{\perp}_xM$ is the normal space of $M$ at $x$ and $\exp^{\perp}$ is the 
normal exponential map of $M$.  Denote by $A$ the shape tensor of $M$ and $R$ 
the curvature tensor of $N$.  
Take $v\in T^{\perp}_xM$ and $X\in T_xM$.  
The strongly $M$-Jacobi field $Y$ along the normal geodesic $\gamma_v$ with 
$Y(0)=X$ (hence $Y'(0)=-A_vX$) is given by 
$$Y(s)=(P_{\gamma_v\vert_{[0,s]}}\circ(D^{co}_{sv}-sD^{si}_{sv}\circ A_v))
(X),$$
where $Y'(0)=\widetilde{\nabla}_vY,\,\,P_{\gamma_v\vert_{[0,s]}}$ is 
the parallel translation along $\gamma_v\vert_{[0,s]}$ and 
$D^{co}_{sv}$ (resp. $D^{si}_{sv}$) is given by 
$$D^{co}_{sv}=\cos(s\sqrt{-1}\sqrt{R(v)})\,\,\,\,\left({\rm resp.}\,\,\,\,
D^{si}_{sv}=\frac{\sin(s\sqrt{-1}\sqrt{R(v)}}
{\sqrt{-1}s\sqrt{R(v)}}\right).$$
All focal radii of $M$ along $\gamma_v$ are obtained as real numbers $s_0$ 
with ${\rm Ker}(D^{co}_{s_0v}-s_0D^{si}_{s_0v}\circ A_v)\not=\{0\}$.  
In general, we ([Koi2]) called a complex number $z_0$ with 
${\rm Ker}(D^{co}_{z_0v}-z_0D^{si}_{z_0v}\circ A_v^{{\bf c}})\not=\{0\}$ 
a {\it complex focal radius of} $M$ {\it along} $\gamma_v$ and call 
${\rm dim}\,{\rm Ker}(D^{co}_{z_0v}-z_0D^{si}_{z_0v}\circ A_v^{{\bf c}})$ the 
{\it multiplicity} of the complex focal radius $z_0$, 
where $A_v^{\bf c}$ is the complexification of $A_v$ and $D^{co}_{z_0v}$ 
(resp. $D^{si}_{z_0v}$) is a ${\bf C}$-linear transformation of 
$(T_xN)^{\bf c}$ defined by 
$$D^{co}_{z_0v}=\cos(z_0\sqrt{-1}\sqrt{R(v)^{\bf c}})\,\,\,\,\left({\rm resp.}
\,\,\,\,D^{si}_{z_0v}=\frac{\sin(z_0\sqrt{-1}\sqrt{R(v)^{\bf c}}}
{z_0\sqrt{-1}\sqrt{R(v)^{\bf c}}}\right).$$
Furthermore, assume that the normal holonomy group of $M$ is trivial.  
Let $v$ be a parallel unit normal vector field of $M$.  
Assume that the number (which may be $0$ and $\infty$) of distinct complex 
focal radii along $\gamma_{v_x}$ is independent of the choice of 
$x\in M$.  Furthermore assume that the number is not equal to $0$.  
Let $\{r_{i,x}\,\vert\,i=1,2,\cdots\}$ 
be the set of all complex focal radii along $\gamma_{v_x}$, where 
$\vert r_{i,x}\vert\,<\,\vert r_{i+1,x}\vert$ or 
"$\vert r_{i,x}\vert=\vert r_{i+1,x}\vert\,\,\&\,\,{\rm Re}\,r_{i,x}
>{\rm Re}\,r_{i+1,x}$" or 
"$\vert r_{i,x}\vert=\vert r_{i+1,x}\vert\,\,\&\,\,
{\rm Re}\,r_{i,x}={\rm Re}\,r_{i+1,x}\,\,\&\,\,
{\rm Im}\,r_{i,x}=-{\rm Im}\,r_{i+1,x}<0$".  
Let $r_i$ ($i=1,2,\cdots$) be complex valued functions on $M$ defined by 
assigning $r_{i,x}$ to each $x\in M$.  
We call these functions $r_i$ ($i=1,2,\cdots$) {\it complex 
focal radius functions for} $v$.  
If, for each parallel unit normal vector field $v$ of $M$, the number 
of distinct complex focal radii along $\gamma_{v_x}$ is independent of 
the choice of $x\in M$, each complex focal radius function for $v$ 
is constant on $M$ and it has constant multiplicity, then 
we call $M$ a {\it complex equifocal submanifold}.  

Next we shall recall the notion of a proper complex isoparametric 
submanifold in a pseudo-Hilbert space.  
Let $M$ be a pseudo-Riemannian Hilbert submanifold in a pseudo-Hilbert space 
$(V,\langle\,\,,\,\,\rangle)$ immersed by $f$.  
See Section 2 of [Koi2] about the definitions of a pseudo-Hilbert space and 
a pseudo-Riemannian Hilbert submanifold.  
Denote by $A$ the shape tensor of $M$ and by $T^{\perp}M$ the normal bundle 
of $M$.  Note that, for $v\in T^{\perp}M$, $A_v$ is not necessarily 
diagonalizable with respect to an orthonormal base and furthermore 
$A_v^{\bf c}$ also is not necessarily diagonalizable with respect to 
a pseudo-orthonormal base.  
We call $M$ a {\it Fredholm pseudo-Riemannian Hilbert submanifold} 
(or simply {\it Fredholm submanifold}) if the following conditions hold:

\vspace{0.2truecm}

(F-i) $M$ is of finite codimension,

(F-ii) There exists an orthogonal time-space decomposition 
$V=V_-\oplus V_+$ such that $(V,\langle\,\,,\,\,\rangle_{V_{\pm}})$ is 
a Hilbert space and that, for each $v\in T^{\perp}M$, $A_v$ is a compact 
operator with respect to $f^{\ast}\langle\,\,,\,\,\rangle_{V_{\pm}}$.  

\vspace{0.2truecm}

\noindent
Since $A_v$ is a compact operator with respect to 
$f^{\ast}\langle\,\,,\,\,\rangle_{V_{\pm}}$, 
the operator ${\rm id}-A_v$ is a Fredholm operator with respect to 
$f^{\ast}\langle\,\,,\,\,\rangle_{V_{\pm}}$ and hence the normal exponential 
map $\exp^{\perp}\,:\,T^{\perp}M\to V$ of 
$M$ is a Fredholm map with respect to the metric of $T^{\perp}M$ naturally 
defined from $f^{\ast}\langle\,\,,\,\,\rangle_{V_{\pm}}$ and 
$\langle\,\,,\,\,\rangle_{V_{\pm}}$, where ${\rm id}$ is the identity 
transformation of $TM$.  
The spectrum of the complexification 
$A_v^{\bf c}$ of $A_v$ is described as 
$\{0\}\cup\{\mu_i\,\vert\,i=1,2,\cdots\}$, where "$\vert\mu_i\vert\,>\,
\vert\mu_{i+1}\vert$" or "$\vert\mu_i\vert=\vert\mu_{i+1}\vert\,\,
\&\,\,{\rm Re}\,\mu_i>{\rm Re}\,\mu_{i+1}$" or 
"$\vert\mu_i\vert=\vert\mu_{i+1}\vert\,\,\&\,\,
{\rm Re}\,\mu_i={\rm Re}\,\mu_{i+1}\,\,\&\,\,{\rm Im}\,\mu_i
=-{\rm Im}\,\mu_{i+1}>0$".  We call 
$\mu_i$ the $i$-{\it th} {\it complex principal curvature of direction} 
$v$.  Assume that the normal holonomy group of $M$ is trivial.  
Let $v$ be a parallel normal vector field on $M$.  
Assume that the number (which may be $\infty$) 
of distinct complex principal curvatures of $v_x$ is 
independent of the choice of $x\in M$.  Then we can define functions 
$\mu_i$ ($i=1,2,\cdots$) on $M$ by assigning the $i$-th complex principal 
curvature of direction $v_x$ to each $x\in M$.  We call this function 
$\mu_i$ the $i$-{\it th} {\it complex principal curvature function 
of direction} $v$.  
If $M$ is a Fredholm submanifold with trivial normal holonomy group 
satisfying the following condition (CI), then we call $M$ a 
{\it complex isoparametric submanifold}:

\vspace{0.2truecm}

\noindent
(CI) $\quad$ for each parallel normal vector field $v$, 
the number of distinct complex principal curvatures of direction 
$v_x$ is independent of the choice of $x\in M$ and 
each complex principal curvature function of direction $v$ 
is constant on $M$ and has constant multiplicity.  

\vspace{0.2truecm}

\noindent
Furthermore, if, for each $v\in T^{\perp}M$, the complexified shape operator 
$A_v^{\bf c}$ is diagonalizable with respect to a pseudo-orthonormal base of 
$(T_xM)^{\bf c}$ ($x\,:\,$ the base point of $v$), that is, there exists 
a pseudo-orthonormal base consisting of the eigenvectors of the complexified 
shape operator $A_v^{\bf c}$, 
then we call $M$ a {\it proper complex isoparametric submanifold}.  
Then, for each $x\in M$, there exists a pseudo-orthonormal base 
of $(T_xM)^{\bf c}$ consisting of the common-eigenvectors of 
the complexified shape operators 
$A_v^{\bf c}$'s ($v\in T^{\perp}_xM$) because 
$A_v^{\bf c}$'s commute.  
Let $\{E_i\,\vert\,i\in I\}$ ($I\subset {\bf N}$) be the family of 
subbundles of $(TM)^{\bf c}$ such that, for each $x\in M$, 
$\{E_i(x)\,\vert\,i\in I\}$ is the set of all common-eigenspaces 
of $A_v^{\bf c}$'s ($v\in T_x^{\perp}M$).  
Note that $\displaystyle{(T_xM)^{\bf c}=
\overline{\mathop{\oplus}_{i\in I}E_i(x)}}$ holds.  
There exist smooth sections $\lambda_i$ ($i\in I$) of 
$((T^{\perp}M)^{\bf c})^{\ast}$ such that 
$A_v^{\bf c}=(\lambda_i)_x(v){\rm id}$ on $(E_i)_x$ for each $x\in M$ and 
each $v\in T_x^{\perp}M$.  
We call $\lambda_i$ ($i\in I$) {\it complex principal curvatures of} $M$ and 
call subbundles $E_i$ ($i\in I$) of $(T^{\perp}M)^{\bf c}$ 
{\it complex curvature distributions of} $M$.  Note that 
$(\lambda_i)_x(v)$ is one of the complex principal curvatures of direction 
$v$.  
Set ${\it l}_i:=(\lambda_i)_x^{-1}(1)\,(\subset(T^{\perp}_xM)^{\bf c})$ and 
$R_i^x$ be the complex reflection of order two with respect to ${\it l}_i^x$, 
where $i\in I$.  Denote by $W^x_M$ the group generated by $R_i^x$'s 
($i\in I$), which is independent of the choice of $x\in M$ up to 
isomorphicness.  We call ${\it l}_i^x$'s {\it complex focal hyperplanes of} 
$(M,x)$.  

\section{Proofs of Theorems A and B} 
In this section, we shall prove Theorems A and B.  
First we prepare a lemma to prove Theorem A.  
Let $M$ be a proper complex equifocal submanifold in a symmetric space $G/K$ 
of non-compact type and $\widetilde M$ a component of 
$(\pi\circ\phi)^{-1}(M)$.  
Without loss of generality, we may assume that $eK$ belongs to $M$, where 
$e$ is the identity element of $G$.  
Hence we may assume that the constant path $\hat 0$ at the zero element $0$ of 
$\mathfrak g$ belongs to $\widetilde M$.  
Fix a unit normal vector $v$ of $M$ at $eK$.  Set $\mathfrak p:=T_{eK}(G/K)$ 
and $\mathfrak b:=T^{\perp}_{eK}M$.  Let $\mathfrak a$ be a maximal abelian 
subspace of $\mathfrak p\,(\subset\mathfrak g)$ containing $\mathfrak b$.  
and $\mathfrak p=\mathfrak a+\sum\limits_{\alpha\in\triangle_+}
\mathfrak p_{\alpha}$ be the root space decomposition with respect to 
$\mathfrak a$, that is, $\mathfrak p_{\alpha}:=\{X\in\mathfrak p\,\vert\,
{\rm ad}(a)^2(X)=\alpha(a)^2X\,\,(\forall\,\,a\in\mathfrak a)\}$ and 
$\triangle_+$ is the positive root system of the root system 
$\triangle:=\{\alpha\in{\mathfrak a}^{\ast}\setminus\{0\}\,\vert\,
\mathfrak p_{\alpha}\not=\{0\}\}$ under a lexicographic ordering of 
${\mathfrak a}^{\ast}$.  
Let $\triangle_{\mathfrak b}:=\{\alpha\vert_{\mathfrak b}\,\vert\,\alpha\in
\triangle\,\,{\rm s.t.}\,\,\alpha\vert_{\mathfrak b}\not=0\}$ and 
$\mathfrak p=\mathfrak z_{\mathfrak p}(\mathfrak b)
+\sum_{\beta\in(\triangle_{\mathfrak b})_+}
\mathfrak p_{\beta}$ be the root space decomposition with respect to 
$\mathfrak b$, where $\mathfrak z_{\mathfrak p}(\mathfrak b)$ is 
the centralizer of $\mathfrak b$ in $\mathfrak p$, 
$\mathfrak p_{\beta}=\sum\limits_{\alpha\in\triangle_+\,{\rm s.t.}\,
\alpha\vert_{\mathfrak b}=\pm\beta}\mathfrak p_{\alpha}$ and 
$(\triangle_{\mathfrak b})_+$ is the positive root system of the root system 
$\triangle_{\mathfrak b}$ under a lexicographic ordering of 
$\mathfrak b^{\ast}$.  
For convenience, we denote $\mathfrak z_{\mathfrak p}(\mathfrak b)$ 
by $\mathfrak p_0$.  
Denote by $A$ (resp. $\widetilde A$) the shape tensor of $M$ (resp. 
$\widetilde M$).  Also, denote by $R$ the curvature tensor of $G/K$.  
Let $m_A:=\displaystyle{\mathop{\max}_{v\in\mathfrak b\setminus\{0\}}
\sharp{\rm Spec}\,A_v}$ and 
$m_R:=\displaystyle{\mathop{\max}_{v\in\mathfrak b\setminus\{0\}}
\sharp{\rm Spec}\,R(v)}$, where $\sharp(\cdot)$ is the cardinal number of 
$(\cdot)$.  Note that $m_R=\sharp(\triangle_{\mathfrak b})_+$.  
Let $U:=\{v\in\mathfrak b\setminus\{0\}\,\vert\,
\sharp{\rm Spec}A_v=m_A,\,\,\sharp{\rm Spec}\,R(v)=m_R\}$, 
which is an open dense subset of $\mathfrak b\setminus\{0\}$.  
Fix $v\in U$.  Note that 
${\rm Spec}\,R(v)=\{-\beta(v)^2\,\vert\,\beta\in(\triangle_{\mathfrak b})_+
\}$.  Since $v\in U$, $\beta(v)^2$'s 
($\beta\in(\triangle_{\mathfrak b})_+$) are mutually distinct.  
Let ${\rm Spec}A_v=\{\lambda^v_1,\cdots,\lambda^v_{m_A}\}$ 
($\lambda^v_1>\cdots>\lambda^v_{m_A}$).  
Set 
$$\begin{array}{l}
\displaystyle{I_0^v:=\{i\,\vert\,\mathfrak p_0\cap{\rm Ker}
(A_v-\lambda_i^v{\rm id})\not=\{0\}\},}\\
\displaystyle{I_{\beta}^v:=\{i\,\vert\,\mathfrak p_{\beta}\cap{\rm Ker}
(A_v-\lambda^v_i{\rm id})\not=\{0\}\},}\\
\displaystyle{(I_{\beta}^v)^+:=\{i\,\vert\,\mathfrak p_{\beta}\cap{\rm Ker}
(A_v-\lambda^v_i{\rm id})\not=\{0\},\,\,\vert\lambda^v_i\vert\,>\,
\vert\beta(v)\vert\},}\\
\displaystyle{(I_{\beta}^v)^-:=\{i\,\vert\,\mathfrak p_{\beta}\cap{\rm Ker}
(A_v-\lambda^v_i{\rm id})\not=\{0\},\,\,
\vert\lambda^v_i\vert\,<\,\vert\beta(v)\vert\},}\\
\displaystyle{(I_{\beta}^v)^0:=\{i\,\vert\,\mathfrak p_{\beta}\cap{\rm Ker}
(A_v-\lambda^v_i{\rm id})\not=\{0\},\,\,
\vert\lambda^v_i\vert=\vert\beta(v)\vert\}.}
\end{array}$$
Let $F$ be the sum of all 
complex focal hyperplanes of $(\widetilde M,\hat 0)$.  
Denote by ${\rm pr}_{\bf R}$ the natural projection of $\mathfrak b^{\bf c}$ 
onto $\mathfrak b$ and set $F_{\bf R}:={\rm pr}_{\bf R}(F)$.  
Then we have the following facts.  

\vspace{0.3truecm}

\noindent
{\bf Lemma 3.1.} {\sl Assume that $M$ is curvature-adapted and that $G/K$ is 
not a hyperbolic space.  Then the set 
$(I_{\beta}^v)^0$ is empty and 
the spectrum of 
$A_v\vert_{\sum_{\beta\in(\triangle_{\mathfrak b})_+}\mathfrak p_{\beta}}$ 
is equal to 
$$\begin{array}{l}
\displaystyle{\{\frac{\beta(v)}{\tanh\beta(Z)}\,\vert\,\beta\in
(\triangle_{\mathfrak b})_+\,\,{\rm s.t.}\,\,(I_{\beta}^v)^+\not=\emptyset\}}\\
\displaystyle{\cup\{\beta(v)\tanh\beta(Z)\,\vert\,\beta\in
(\triangle_{\mathfrak b})_+\,\,{\rm s.t.}\,\,(I_{\beta}^v)^-\not=\emptyset\}}
\end{array}$$
for some $Z\in\mathfrak b$.}

\vspace{0.3truecm}

\noindent
{\it Proof.} From $v\in U$, we have $\beta(v)\not=0$ for any 
$\beta\in(\triangle_{\mathfrak b})_+$.  Hence, since $M$ is 
curvature-adapted and proper complex equifocal, it follows from Theorem 1 of 
[Koi3] that $(I_{\beta}^v)^0=\emptyset$.  
Set $c^+_{\beta,i,v}:=\frac{\beta(v)}{\lambda_i^v}$ 
$(i\in(I_{\beta}^v)^+\,(\beta\in(\triangle_{\mathfrak b})_+))$ 
and $c^-_{\beta,i,v}:=\frac{\lambda_i^v}{\beta(v)}$ 
$(i\in(I_{\beta}^v)^-\,(\beta\in(\triangle_{\mathfrak b})_+))$.  
According to the proof of Theorems B and C in [Koi11], we have 
$$\begin{array}{l}
\displaystyle{F=\left(
\mathop{\cup}_{\beta\in(\triangle_{\mathfrak b})_+}
\mathop{\cup}_{(i,j)\in(I^v_{\beta})^+\times{\bf Z}}
(\beta^{\bf c})^{-1}({\rm arctanh}c_{\beta,i,v}^++j\pi\sqrt{-1})\right)}\\
\hspace{1truecm}\displaystyle{\cup\left(
\mathop{\cup}_{\beta\in(\triangle_{\mathfrak b})_+}
\mathop{\cup}_{(i,j)\in(I^v_{\beta})^-\times{\bf Z}}
(\beta^{\bf c})^{-1}({\rm arctanh}c_{\beta,i,v}^-+(j+\frac12)\pi\sqrt{-1})
\right)}
\end{array}
\leqno{(3.1)}$$
and 
$$\begin{array}{l}
\displaystyle{F_{\bf R}=\left(
\mathop{\cup}_{\beta\in(\triangle_{\mathfrak b})_+}
\mathop{\cup}_{i\in(I^v_{\beta})^+}
\beta^{-1}({\rm arctanh}c_{\beta,i,v}^+)\right)}\\
\hspace{1truecm}\displaystyle{\cup\left(
\mathop{\cup}_{\beta\in(\triangle_{\mathfrak b})_+}
\mathop{\cup}_{i\in(I^v_{\beta})^-}
\beta^{-1}({\rm arctanh}c_{\beta,i,v}^-)\right).}
\end{array}$$
Also, since $G/K$ is not a hyperbolic space, the intersection of all 
the hyperplanes constructing $F_{\bf R}$ is non-empty.  
Take an element $Z$ of the intersection.  
Then we have 
$$\lambda_i^v=\left\{
\begin{array}{ll}
\displaystyle{\frac{\beta(v)}{\tanh\beta(Z)}}&(i\in(I_{\beta}^v)^+)\\
\displaystyle{\beta(v)\tanh\beta(Z)}&(i\in(I_{\beta}^v)^-).
\end{array}\right.$$
Hence we have 
$$\begin{array}{l}
\displaystyle{{\rm Spec}
\left(A_v\vert_{\sum_{\beta\in(\triangle_{\mathfrak b})_+}
\mathfrak p_{\beta}}\right)=\{\frac{\beta(v)}{\tanh\beta(Z)}\,\vert\,\beta\in
(\triangle_{\mathfrak b})_+\,\,{\rm s.t.}\,\,(I_{\beta}^v)^+\not=\emptyset\}}\\
\hspace{4.5truecm}
\displaystyle{\cup\{\beta(v)\tanh\beta(Z)\,\vert\,\beta\in
(\triangle_{\mathfrak b})_+\,\,{\rm s.t.}\,\,(I_{\beta}^v)^-\not=\emptyset\}.}
\end{array}$$
\begin{flushright}q.e.d.\end{flushright}

\vspace{0.3truecm}

By using this lemma, we prove Theorem A.  

\vspace{0.3truecm}

\noindent
{\it Proof of Theorem A.} 
Let $M$ be as in the statement of Theorem A.  
In the case where $G/K$ is a hyperbolic space, $M$ is an isoparametric 
hypersurface (other than a horosphere) in the space and hence it occurs as 
a principal orbit of a Hermann type action.  So we suffice to show 
the statement in the case where $G/K$ is not a hyperbolic space.  
Let $Z$ be as in the proof of Lemma 3.1 and 
$\widetilde Z$ be the parallel normal vector field of $M$ with 
$\widetilde Z_{eK}=Z$.  Denote by $\eta_{t\widetilde Z}$ ($0\leq t\leq1$) the 
end-point map for $t\widetilde Z$ (i.e., $\eta_{t\widetilde Z}(x)
=\exp^{\perp}(t\widetilde Z_x)\,\,(x\in M)$), where $\exp^{\perp}$ is the 
normal exponential map of $M$.  Set $M_t:=\eta_{t\widetilde Z}(M)$, which is 
a parallel submanifold or a focal submanifold of $M$.  In particular, $M_1$ is 
a focal submanifold of $M$.  Let $\widetilde v$ be the parallel 
tangent vector field on the flat section $\Sigma_{eK}$ with 
$\widetilde v_{eK}=v(\in U)$.  Note that 
$\widetilde v_{\eta_{t\widetilde Z}(eK)}$ is a normal vector of $M_t$.  
Denote by $A^t$ the shape tensors of $M_t$ ($0<t\leq1$).  
It is clear that such that $M_t$'s ($1-\varepsilon<t<1$) are parallel 
submanifolds of $M$ for a sufficiently small positive number $\varepsilon$.  
According to the proof of Theorems B and C of [Koi11], there 
exists a complex linear function $\phi_i$ on $\mathfrak b^{\bf c}$ with 
$\phi_i(v)=\lambda_i^v$ and $\phi_i^{-1}(1)\subset F$ for each $i\in I_0^v$ 
with $\lambda_i^v\not=0$.  
Fix $i_0\in I_0^v$.  According to $(3.1)$, $\phi_{i_0}^{-1}(1)$ coincides with 
one of $(\beta^{\bf c})^{-1}({\rm arctanh}c_{\beta,i,v}^++j\pi\sqrt{-1})$'s 
($\beta\in(\triangle_{\mathfrak b})_+,\,(i,j)\in(I_{\beta}^v)^+\times{\bf Z}$) 
and 
$(\beta^{\bf c})^{-1}({\rm arctanh}c_{\beta,i,v}^-+(j+\frac12)\pi\sqrt{-1})$'s 
($\beta\in(\triangle_{\mathfrak b})_+,\,(i,j)\in
(I_{\beta}^v)^-\times{\bf Z}$).  Since 
$(\beta^{\bf c})^{-1}({\rm arctanh}c_{\beta,i,v}^++j\pi\sqrt{-1})=
\beta^{-1}({\rm arctanh}c_{\beta,i,v}^+)+
(\beta^{\bf c}\vert_{\sqrt{-1}\mathfrak b})^{-1}(j\pi\sqrt{-1})$ and 
$(\beta^{\bf c})^{-1}({\rm arctanh}c_{\beta,i,v}^-+(j+\frac12)\pi\sqrt{-1})=
\beta^{-1}({\rm arctanh}c_{\beta,i,v}^-)+
(\beta^{\bf c}\vert_{\sqrt{-1}\mathfrak b})^{-1}((j+\frac12)\pi\sqrt{-1})$, 
we have 
$\phi_{i_0}^{-1}(1)=\beta_1^{-1}({\rm arctanh}c_{\beta_1,i_1,v}^+)+
(\beta_1^{\bf c}\vert_{\sqrt{-1}\mathfrak b})^{-1}(0)$ for some 
$\beta_1\in(\triangle_{\mathfrak b})_+$ and $i_1\in(I_{\beta_1}^v)^+$.  
Hence we have $Z\in\phi_{i_0}^{-1}(1)$, which implies that 
$$\eta_{\widetilde Z\ast}(\mathfrak p_0\cap{\rm Ker}(A_v-\lambda_{i_0}^v
{\rm id}))=\{0\}.\leqno{(3.2)}$$
Easily we can show 
$$A^1_{\widetilde v_{\eta_{\widetilde Z}(eK)}}\vert_
{(\eta_{\widetilde Z})_{\ast}(\mathfrak p_0\cap{\rm Ker}\,A_v)}=0.\leqno{(3.3)}
$$
By discussing delicately in terms of Lemma 3.1, we can show 
$$\begin{array}{l}
\displaystyle{{\rm Spec}\,A^t_{\widetilde v_{\eta_{t\widetilde Z}(eK)}}
\vert_{\sum_{\beta\in\overline{\triangle}_+}(\eta_{t\widetilde Z})_{\ast}
(\mathfrak p_{\beta})}}\\
\displaystyle{=\{\frac{\beta(v)}{\tanh((1-t)\beta(Z))}\,\vert\,\beta\in
(\triangle_{\mathfrak b})_+\,\,{\rm s.t.}\,\,(I_{\beta}^v)^+\not=\emptyset\}}\\
\hspace{0.6truecm}
\displaystyle{\cup\{\beta(v)\tanh((1-t)\beta(Z))\,\vert\,\beta\in
(\triangle_{\mathfrak b})_+\,\,{\rm s.t.}\,\,(I_{\beta}^v)^-\not=\emptyset\}}
\end{array}
\leqno{(3.4)}$$
for each $t\in(1-\varepsilon,1)$.  Denote by 
${\rm Gr}_m(G/K)$ the Grassmann bundle of $G/K$ consisting of $m$-dimensional 
subspaces of the tangent spaces, where $m:={\rm dim}\,M_1$.  
Define $D_t$ ($1-\varepsilon<t<1$) by 
$$\begin{array}{l}
\displaystyle{D_t:=(\eta_{t\widetilde Z})_{\ast}
(\mathfrak p_0\cap{\rm Ker}A_v)}\\
\hspace{1.2truecm}
\displaystyle{+\sum_{\beta\in(\triangle_{\mathfrak b})_+\,\,{\rm s.t.}\,\,
(I_{\beta}^v)^-\not=\emptyset}\left(\mathfrak p_{\beta}\cap{\rm Ker}\left(
A^t_{\widetilde v_{\eta_{t\widetilde Z}(eK)}}-\beta(v)\tanh((1-t)\beta(Z))
{\rm id}\right)\right).}
\end{array}$$
From $(3.2),(3.3)$ and $(3.4)$, we can show $D_t\in{\rm Gr}_m(G/K)$, 
$\lim\limits_{t\to1-0}D_t=T_{\eta_{\widetilde Z}(eK)}M_1$ 
(in ${\rm Gr}_m(G/K)$) and 
$${\rm Spec}\,A^1_{\widetilde v_{\eta_{\widetilde Z}(eK)}}=
\{0\}\cup\{\lim_{t\to1-0}\beta(v)\tanh((1-t)\beta(Z))\,\vert\,\beta\in
(\triangle_{\mathfrak b})_+\,\,{\rm s.t.}\,\,(I_{\beta}^v)^-\not=\emptyset\}=\{0\}.
$$
Thus we have $A^1_{\widetilde v_{\eta_{\widetilde Z}(eK)}}=0$.  
Since this relation holds for any $v\in U$ and $U$ is open and dense in 
$\mathfrak b(=T^{\perp}_{eK}M)$, 
$A^1_{\widetilde v_{\eta_{\widetilde Z}(eK)}}=0$ holds for any 
$v\in T^{\perp}_{eK}M$.  Set 
$L:=\eta_{\widetilde Z}^{-1}(\eta_{\widetilde Z}(eK))$.  Take any $x\in L$.  
Similarly we can show $A^1_{\widetilde v_{\eta_{\widetilde Z}(x)}}=0$ for any 
$v\in T^{\perp}_xM$, where $\widetilde v$ is the parallel tangent vector field 
on the section $\Sigma_x$ of $M$ through $x$ with $\widetilde v_x=v$.  
It is easy to show that $\eta_{\widetilde Z}(x)=\eta_{\widetilde Z}(eK)$ and 
$\sum\limits_{x\in L}\{\widetilde v_{\eta_{\widetilde Z}(x)}\,\vert\,
v\in T^{\perp}_xM\}=T^{\perp}_{\eta_{\widetilde Z}(eK)}M_1$.  Hence we 
see that $A^1$ vanishes at $\eta_{\widetilde Z}(eK)$.  
Similarly we can show that 
$A^1$ vanishes at any point of $M_1$ other than $\eta_{\widetilde Z}(eK)$.  
Therefore $M_1$ is totally geodesic in $G/K$.  
By the way, since $M$ is irreducible, proper complex equifocal and 
${\rm codim}\,M\geq2$, 
it follows from the homogeneity theorem of [Koi10] that $M$ is homogeneous.  
Hence it follows from Theorem A of [Koi6] that $M$ is 
a principal orbit of a complex hyperpolar action on $G/K$.  
Furthermore, since this action admits 
a totally geodesic singular orbit $M_1$ and it is of cohomogeneity greater 
than one, it follows from Theorem C and Remark 1.1 of [Koi6] that it is orbit 
equivalent to a Hermann type action.  Thus we see that $M$ is 
a principal orbit of a Hermann type action.  \hspace{11.8truecm}q.e.d.

\vspace{0.4truecm}

\centerline{
\unitlength 0.1in
\begin{picture}( 66.6000, 26.9900)(-13.4000,-38.0900)
%
\special{pn 8}%
\special{ar 3690 1840 400 200  0.1721908 0.2001811}%
%
\special{pn 8}%
\special{ar 2790 4030 930 2210  5.0709435 6.1828940}%
%
\special{pn 8}%
\special{pa 3000 2610}%
\special{pa 2670 3430}%
\special{fp}%
%
\special{pn 8}%
\special{pa 2670 3430}%
\special{pa 5010 2790}%
\special{fp}%
%
\special{pn 8}%
\special{pa 5010 2790}%
\special{pa 5320 2020}%
\special{fp}%
\special{pa 5320 2020}%
\special{pa 5320 2020}%
\special{fp}%
%
\special{pn 8}%
\special{pa 5320 2030}%
\special{pa 4280 2290}%
\special{fp}%
%
\special{pn 8}%
\special{pa 3010 2610}%
\special{pa 3430 2510}%
\special{fp}%
%
\special{pn 8}%
\special{pa 3490 2490}%
\special{pa 4220 2310}%
\special{fp}%
%
\special{pn 8}%
\special{pa 3490 3200}%
\special{pa 4430 2250}%
\special{fp}%
%
\special{pn 20}%
\special{sh 1}%
\special{ar 3750 2940 10 10 0  6.28318530717959E+0000}%
\special{sh 1}%
\special{ar 3750 2940 10 10 0  6.28318530717959E+0000}%
%
\special{pn 20}%
\special{sh 1}%
\special{ar 4360 2620 10 10 0  6.28318530717959E+0000}%
\special{sh 1}%
\special{ar 4360 2620 10 10 0  6.28318530717959E+0000}%
%
\special{pn 8}%
\special{pa 4500 3460}%
\special{pa 4498 3492}%
\special{pa 4486 3522}%
\special{pa 4470 3550}%
\special{pa 4448 3572}%
\special{pa 4424 3592}%
\special{pa 4396 3610}%
\special{pa 4368 3626}%
\special{pa 4340 3640}%
\special{pa 4310 3652}%
\special{pa 4280 3662}%
\special{pa 4248 3670}%
\special{pa 4218 3678}%
\special{pa 4186 3684}%
\special{pa 4154 3690}%
\special{pa 4122 3692}%
\special{pa 4090 3694}%
\special{pa 4058 3696}%
\special{pa 4026 3694}%
\special{pa 3994 3692}%
\special{pa 3962 3690}%
\special{pa 3932 3684}%
\special{pa 3900 3678}%
\special{pa 3868 3670}%
\special{pa 3838 3660}%
\special{pa 3810 3648}%
\special{pa 3780 3632}%
\special{pa 3756 3612}%
\special{pa 3732 3590}%
\special{pa 3714 3564}%
\special{pa 3704 3534}%
\special{pa 3704 3502}%
\special{pa 3714 3470}%
\special{pa 3730 3444}%
\special{pa 3750 3420}%
\special{pa 3774 3398}%
\special{pa 3802 3380}%
\special{pa 3828 3364}%
\special{pa 3858 3350}%
\special{pa 3886 3338}%
\special{pa 3918 3326}%
\special{pa 3948 3318}%
\special{pa 3978 3310}%
\special{pa 4010 3302}%
\special{pa 4042 3298}%
\special{pa 4074 3294}%
\special{pa 4106 3292}%
\special{pa 4138 3290}%
\special{pa 4170 3290}%
\special{pa 4202 3292}%
\special{pa 4234 3296}%
\special{pa 4266 3300}%
\special{pa 4296 3306}%
\special{pa 4328 3314}%
\special{pa 4358 3324}%
\special{pa 4388 3336}%
\special{pa 4416 3350}%
\special{pa 4442 3368}%
\special{pa 4466 3390}%
\special{pa 4486 3416}%
\special{pa 4498 3444}%
\special{pa 4500 3460}%
\special{sp}%
%
\special{pn 8}%
\special{pa 4340 2620}%
\special{pa 4344 2652}%
\special{pa 4338 2684}%
\special{pa 4326 2714}%
\special{pa 4308 2740}%
\special{pa 4286 2764}%
\special{pa 4262 2786}%
\special{pa 4238 2806}%
\special{pa 4212 2824}%
\special{pa 4186 2842}%
\special{pa 4156 2856}%
\special{pa 4128 2870}%
\special{pa 4098 2882}%
\special{pa 4068 2894}%
\special{pa 4038 2904}%
\special{pa 4006 2912}%
\special{pa 3976 2920}%
\special{pa 3944 2926}%
\special{pa 3912 2930}%
\special{pa 3882 2934}%
\special{pa 3850 2936}%
\special{pa 3818 2938}%
\special{pa 3786 2936}%
\special{pa 3754 2932}%
\special{pa 3722 2928}%
\special{pa 3690 2920}%
\special{pa 3660 2910}%
\special{pa 3632 2896}%
\special{pa 3606 2878}%
\special{pa 3584 2854}%
\special{pa 3568 2826}%
\special{pa 3562 2794}%
\special{pa 3566 2762}%
\special{pa 3578 2734}%
\special{pa 3594 2706}%
\special{pa 3614 2680}%
\special{pa 3638 2658}%
\special{pa 3662 2638}%
\special{pa 3688 2618}%
\special{pa 3716 2602}%
\special{pa 3744 2586}%
\special{pa 3772 2572}%
\special{pa 3800 2558}%
\special{pa 3830 2548}%
\special{pa 3862 2538}%
\special{pa 3892 2528}%
\special{pa 3922 2520}%
\special{pa 3954 2514}%
\special{pa 3986 2510}%
\special{pa 4018 2506}%
\special{pa 4050 2504}%
\special{pa 4082 2502}%
\special{pa 4114 2502}%
\special{pa 4146 2504}%
\special{pa 4178 2510}%
\special{pa 4208 2516}%
\special{pa 4240 2526}%
\special{pa 4268 2540}%
\special{pa 4296 2556}%
\special{pa 4318 2578}%
\special{pa 4334 2606}%
\special{pa 4340 2620}%
\special{sp}%
%
\special{pn 8}%
\special{pa 4034 1844}%
\special{pa 4044 1874}%
\special{pa 4044 1906}%
\special{pa 4038 1938}%
\special{pa 4026 1968}%
\special{pa 4012 1996}%
\special{pa 3992 2022}%
\special{pa 3972 2046}%
\special{pa 3950 2070}%
\special{pa 3928 2092}%
\special{pa 3902 2112}%
\special{pa 3878 2132}%
\special{pa 3852 2150}%
\special{pa 3824 2168}%
\special{pa 3796 2184}%
\special{pa 3768 2198}%
\special{pa 3738 2212}%
\special{pa 3710 2224}%
\special{pa 3680 2236}%
\special{pa 3648 2246}%
\special{pa 3618 2254}%
\special{pa 3586 2262}%
\special{pa 3556 2268}%
\special{pa 3524 2270}%
\special{pa 3492 2272}%
\special{pa 3460 2270}%
\special{pa 3428 2266}%
\special{pa 3396 2258}%
\special{pa 3368 2246}%
\special{pa 3342 2228}%
\special{pa 3320 2204}%
\special{pa 3308 2174}%
\special{pa 3306 2142}%
\special{pa 3312 2110}%
\special{pa 3322 2080}%
\special{pa 3338 2052}%
\special{pa 3354 2024}%
\special{pa 3374 2000}%
\special{pa 3396 1976}%
\special{pa 3420 1954}%
\special{pa 3444 1932}%
\special{pa 3470 1914}%
\special{pa 3496 1894}%
\special{pa 3522 1876}%
\special{pa 3548 1860}%
\special{pa 3578 1846}%
\special{pa 3606 1832}%
\special{pa 3636 1818}%
\special{pa 3666 1808}%
\special{pa 3696 1798}%
\special{pa 3726 1790}%
\special{pa 3758 1782}%
\special{pa 3790 1776}%
\special{pa 3820 1772}%
\special{pa 3852 1768}%
\special{pa 3884 1768}%
\special{pa 3916 1772}%
\special{pa 3948 1778}%
\special{pa 3978 1790}%
\special{pa 4006 1808}%
\special{pa 4026 1830}%
\special{pa 4034 1844}%
\special{sp}%
%
\special{pn 8}%
\special{ar 3240 4250 1320 3060  5.1600065 6.1237801}%
%
\special{pn 8}%
\special{ar 2390 3790 1910 2900  5.4241963 6.2783705}%
%
\special{pn 8}%
\special{ar 2360 3770 1560 2500  5.3587035 6.2831853}%
\special{ar 2360 3770 1560 2500  0.0000000 0.0020906}%
%
\special{pn 8}%
\special{pa 4150 2670}%
\special{pa 4148 2702}%
\special{pa 4130 2730}%
\special{pa 4108 2750}%
\special{pa 4080 2768}%
\special{pa 4052 2784}%
\special{pa 4022 2794}%
\special{pa 3992 2804}%
\special{pa 3960 2810}%
\special{pa 3928 2814}%
\special{pa 3896 2814}%
\special{pa 3864 2812}%
\special{pa 3834 2802}%
\special{pa 3806 2786}%
\special{pa 3790 2758}%
\special{pa 3796 2728}%
\special{pa 3812 2702}%
\special{pa 3836 2680}%
\special{pa 3862 2662}%
\special{pa 3892 2648}%
\special{pa 3922 2636}%
\special{pa 3952 2628}%
\special{pa 3984 2622}%
\special{pa 4016 2618}%
\special{pa 4048 2618}%
\special{pa 4080 2620}%
\special{pa 4110 2630}%
\special{pa 4136 2648}%
\special{pa 4150 2670}%
\special{sp}%
%
\special{pn 8}%
\special{pa 4280 3460}%
\special{pa 4274 3492}%
\special{pa 4254 3516}%
\special{pa 4228 3534}%
\special{pa 4198 3548}%
\special{pa 4170 3560}%
\special{pa 4138 3568}%
\special{pa 4106 3574}%
\special{pa 4074 3574}%
\special{pa 4042 3576}%
\special{pa 4010 3572}%
\special{pa 3980 3566}%
\special{pa 3950 3552}%
\special{pa 3926 3532}%
\special{pa 3914 3504}%
\special{pa 3922 3474}%
\special{pa 3942 3448}%
\special{pa 3968 3430}%
\special{pa 3998 3418}%
\special{pa 4028 3406}%
\special{pa 4058 3398}%
\special{pa 4090 3394}%
\special{pa 4122 3392}%
\special{pa 4154 3392}%
\special{pa 4186 3396}%
\special{pa 4216 3402}%
\special{pa 4246 3416}%
\special{pa 4270 3436}%
\special{pa 4280 3460}%
\special{sp}%
%
\special{pn 8}%
\special{pa 3840 1920}%
\special{pa 3844 1952}%
\special{pa 3834 1982}%
\special{pa 3814 2008}%
\special{pa 3792 2030}%
\special{pa 3768 2050}%
\special{pa 3740 2068}%
\special{pa 3712 2084}%
\special{pa 3682 2096}%
\special{pa 3652 2106}%
\special{pa 3622 2112}%
\special{pa 3590 2116}%
\special{pa 3558 2112}%
\special{pa 3528 2102}%
\special{pa 3506 2078}%
\special{pa 3504 2048}%
\special{pa 3516 2018}%
\special{pa 3536 1992}%
\special{pa 3558 1970}%
\special{pa 3584 1950}%
\special{pa 3610 1934}%
\special{pa 3640 1918}%
\special{pa 3668 1906}%
\special{pa 3698 1896}%
\special{pa 3730 1890}%
\special{pa 3762 1886}%
\special{pa 3794 1890}%
\special{pa 3824 1902}%
\special{pa 3840 1920}%
\special{sp}%
%
\special{pn 8}%
\special{ar 2280 3930 1840 2970  5.4215498 6.2389105}%
%
\special{pn 20}%
\special{sh 1}%
\special{ar 4100 3480 10 10 0  6.28318530717959E+0000}%
\special{sh 1}%
\special{ar 4100 3480 10 10 0  6.28318530717959E+0000}%
%
\special{pn 20}%
\special{sh 1}%
\special{ar 3680 2010 10 10 0  6.28318530717959E+0000}%
\special{sh 1}%
\special{ar 3680 2010 10 10 0  6.28318530717959E+0000}%
%
\special{pn 20}%
\special{sh 1}%
\special{ar 3970 2720 10 10 0  6.28318530717959E+0000}%
\special{sh 1}%
\special{ar 3970 2720 10 10 0  6.28318530717959E+0000}%
%
\special{pn 8}%
\special{pa 3600 3430}%
\special{pa 3750 2960}%
\special{dt 0.045}%
\special{sh 1}%
\special{pa 3750 2960}%
\special{pa 3712 3018}%
\special{pa 3734 3012}%
\special{pa 3750 3030}%
\special{pa 3750 2960}%
\special{fp}%
%
\special{pn 8}%
\special{pa 2850 3000}%
\special{pa 5170 2410}%
\special{fp}%
\special{pa 4260 2760}%
\special{pa 4260 2760}%
\special{fp}%
%
\special{pn 8}%
\special{pa 4960 2330}%
\special{pa 4860 2480}%
\special{dt 0.045}%
\special{sh 1}%
\special{pa 4860 2480}%
\special{pa 4914 2436}%
\special{pa 4890 2436}%
\special{pa 4880 2414}%
\special{pa 4860 2480}%
\special{fp}%
%
\special{pn 8}%
\special{pa 3440 3060}%
\special{pa 3560 3110}%
\special{dt 0.045}%
\special{sh 1}%
\special{pa 3560 3110}%
\special{pa 3506 3066}%
\special{pa 3512 3090}%
\special{pa 3492 3104}%
\special{pa 3560 3110}%
\special{fp}%
%
\special{pn 13}%
\special{pa 3910 2450}%
\special{pa 4050 2980}%
\special{fp}%
%
\special{pn 13}%
\special{pa 4090 2390}%
\special{pa 4230 2920}%
\special{fp}%
%
\special{pn 13}%
\special{pa 4300 2330}%
\special{pa 4440 2860}%
\special{fp}%
%
\special{pn 8}%
\special{pa 4640 2000}%
\special{pa 4110 2440}%
\special{dt 0.045}%
\special{sh 1}%
\special{pa 4110 2440}%
\special{pa 4174 2414}%
\special{pa 4152 2406}%
\special{pa 4150 2382}%
\special{pa 4110 2440}%
\special{fp}%
%
\special{pn 8}%
\special{pa 4910 1980}%
\special{pa 4340 2460}%
\special{dt 0.045}%
\special{sh 1}%
\special{pa 4340 2460}%
\special{pa 4404 2432}%
\special{pa 4382 2426}%
\special{pa 4378 2402}%
\special{pa 4340 2460}%
\special{fp}%
%
\special{pn 8}%
\special{pa 4560 2470}%
\special{pa 4380 2600}%
\special{dt 0.045}%
\special{sh 1}%
\special{pa 4380 2600}%
\special{pa 4446 2578}%
\special{pa 4424 2570}%
\special{pa 4422 2546}%
\special{pa 4380 2600}%
\special{fp}%
\put(37.3000,-11.5000){\makebox(0,0)[rt]{$M_1$}}%
\put(40.1000,-11.1000){\makebox(0,0)[rt]{$M_t$}}%
\put(48.2000,-11.6000){\makebox(0,0)[rt]{$M$}}%
\put(47.5000,-23.3000){\makebox(0,0)[rt]{$eK$}}%
\put(51.4000,-21.6000){\makebox(0,0)[rt]{$\Sigma_{eK}$}}%
\put(36.4000,-34.9000){\makebox(0,0)[rt]{$x$}}%
\put(33.8000,-29.8000){\makebox(0,0)[rt]{$\Sigma_x$}}%
\put(51.3000,-18.0000){\makebox(0,0)[rt]{$D_0$}}%
\put(48.2000,-18.2000){\makebox(0,0)[rt]{$D_t$}}%
%
\special{pn 8}%
\special{pa 4470 1740}%
\special{pa 3920 2470}%
\special{dt 0.045}%
\special{sh 1}%
\special{pa 3920 2470}%
\special{pa 3976 2430}%
\special{pa 3952 2428}%
\special{pa 3944 2406}%
\special{pa 3920 2470}%
\special{fp}%
\put(50.0000,-15.4000){\makebox(0,0)[rt]{$T_{\eta_{\widetilde Z}(eK)}M_1$}}%
%
\special{pn 8}%
\special{pa 3030 2280}%
\special{pa 3950 2710}%
\special{dt 0.045}%
\special{sh 1}%
\special{pa 3950 2710}%
\special{pa 3898 2664}%
\special{pa 3902 2688}%
\special{pa 3882 2700}%
\special{pa 3950 2710}%
\special{fp}%
\put(29.8000,-22.0000){\makebox(0,0)[rt]{$\eta_{\widetilde Z}(eK)=\eta_{\widetilde Z}(x)$}}%
%
\special{pn 8}%
\special{pa 4590 1320}%
\special{pa 3860 1700}%
\special{dt 0.045}%
\special{sh 1}%
\special{pa 3860 1700}%
\special{pa 3928 1688}%
\special{pa 3908 1676}%
\special{pa 3910 1652}%
\special{pa 3860 1700}%
\special{fp}%
%
\special{pn 8}%
\special{pa 3820 1300}%
\special{pa 3620 1720}%
\special{dt 0.045}%
\special{sh 1}%
\special{pa 3620 1720}%
\special{pa 3668 1668}%
\special{pa 3644 1672}%
\special{pa 3632 1652}%
\special{pa 3620 1720}%
\special{fp}%
%
\special{pn 8}%
\special{pa 3540 1340}%
\special{pa 3520 1750}%
\special{dt 0.045}%
\special{sh 1}%
\special{pa 3520 1750}%
\special{pa 3544 1684}%
\special{pa 3524 1698}%
\special{pa 3504 1682}%
\special{pa 3520 1750}%
\special{fp}%
%
\special{pn 8}%
\special{pa 2690 2490}%
\special{pa 3620 2660}%
\special{dt 0.045}%
\special{sh 1}%
\special{pa 3620 2660}%
\special{pa 3558 2628}%
\special{pa 3568 2650}%
\special{pa 3552 2668}%
\special{pa 3620 2660}%
\special{fp}%
\put(26.7000,-24.2000){\makebox(0,0)[rt]{$L$}}%
\end{picture}%
\hspace{11truecm}}

\vspace{0.5truecm}

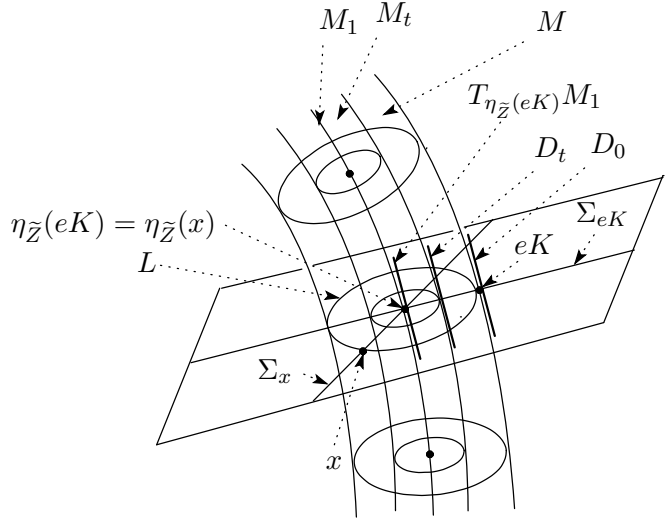
\centerline{{\bf Fig. 1.}}

\vspace{0.5truecm}

Next we shall prove Theorems B and C.  
For its purpose, we prepare the following lemma.  

\vspace{0.5truecm}

\noindent
{\bf Lemma 3.2.} {\sl Let $\psi_1$ and $\psi_2$ be a skew-symmetric 
${\bf C}$-linear transformation and a symmetric ${\bf C}$-linear 
transformation of a (finite dimensional) anti-Kaehlerian space 
$(V,\langle\,\,,\,\,\rangle)$, respectively.  Take $X\in V\setminus\{0\}$.  
Assume that 
$$\left(\cosh((a+bk)\psi_1)-\frac{\sinh((a+bk)\psi_1)}{\psi_1}\circ\psi_2
\right)(X)=0$$
for all $k\in{\bf Z}$, where $a\in{\bf C}$ and $b\in{\bf C}\setminus\{0\}$.  
Let $S$ be the minimal subset of the spectrum of $\psi_1^2$ satisfying 
$\{X,\psi_2(X)\}\subset\displaystyle{\mathop{\oplus}_{\mu\in S}{\rm Ker}
(\psi_1^2-\mu\,{\rm id})}$.  Then we have $b\sqrt{\mu}\in\pi\sqrt{-1}{\bf Z}$ 
for any $\mu\in S$.}

\vspace{0.5truecm}

\noindent
{\it Proof.} Let $X=\sum_{\mu\in S}X_{\mu}$ and $\psi_2(X)=\sum_{\mu\in S}
\psi_2(X)_{\mu}$, where $X_{\mu},\,\psi_2(X)_{\mu}\in{\rm Ker}(\psi_1^2-
\mu\,{\rm id})$.  Thus it follows from the assumption that 
$$\cosh((a+bk)\sqrt{\mu})X_{\mu}-\frac{\sinh((a+bk)\sqrt{\mu})}{\sqrt{\mu}}
\psi_2(X)_{\mu}=0\,\,\,\,(\mu\in S,\,k\in{\bf Z}).$$
From the minimality of $S$, either $X_{\mu}$ or $\psi_2(X)_{\mu}$ is not equal 
to the zero vector.  Hence we have $b\sqrt{\mu}\in \pi\sqrt{-1}{\bf Z}$.  
\hspace{10.2truecm}q.e.d.

\vspace{0.5truecm}

By using Lemma 3.2, we prove Theorem B.  

\vspace{0.5truecm}

\noindent
{\it Proof of Theorem B.} 
Let $M$ be a proper complex equifocal submanifold in 
a symmetric space $G/K$ of non-compact type and $\widetilde M$ 
a component of the lifted submanifold $(\pi\circ\phi)^{-1}(M)$.  
Let $\{\lambda_i\,\vert\,i\in I\}$ be the set of all complex principal 
curvatures of $\widetilde M$ and set 
${\it l}^u_i:=(\lambda_i)_u^{-1}(1)$ ($u\in\widetilde M,\,i\in I$).  Since the 
group generated by the complex reflections $R^u_i$'s ($i\in I$) of order $2$ 
with respect to ${\it l}^u_i$ is discrete, ${\cal L}_u:=\{{\it l}^u_i\,\vert\,i
\in I\}$ is written in the form of the sum 
$\displaystyle{\mathop{\cup}_{j=1}^r{\cal L}^u_j}$ of subfamilies 
${\cal L}^u_j:=\{{\it l}^u_{i(j,k)}\,\vert\,k\in{\bf Z}\}$ ($j=1,\cdots,r$) of 
parallel complex hyperplanes equidistant to one another.  
For each $j\in\{1,\cdots,r\}$, we have $(\lambda_{i(j,k)})_u^{-1}(1)
=(\lambda_{i(j,0)})_u^{-1}(1+kb_j)$ for some $b_j\in{\bf C}$ which is 
independent of the choice of $u\in\widetilde M$.  For simplicity, 
we denote $\lambda_{i(j,0)}$ by $\overline{\lambda}_j$ $(j=1,\cdots,r$).  

\noindent
{\bf (Step I)} First we shall show that, for each parallel normal vector field 
$v$ of $M$, ${\rm Spec}\,R(v_x)$ is independent of the choice of $x\in M$.  
Fix $x\in M$.  Denote by $v^L$ the horizontal lift of $v$ to 
$\widetilde M$.  Since $M$ has flat section, $v^L$ is parallel with respect 
to the normal connection of $\widetilde M$.  
The set $FR_{v_x}$ of all complex focal radii of $M$ 
along $\gamma_{v_x}$ is given by 
$$FR_{v_x}=\{z\in{\bf C}\,\vert\,{\rm Ker}(D^{co}_{zv_x}
-D^{si}_{zv_x}\circ A_{zv_x}^{\bf c})\not=\{0\}\}.$$
For simplicity, set $Q_{v_x}(z):=D^{co}_{zv_x}-D^{si}_{zv_x}\circ 
A^{\bf c}_{zv_x}$.  
Take $u\in\widetilde M\cap(\pi\circ\phi)^{-1}(x)$.  
The set $FR_{v^L_u}$ of all complex focal radii of $\widetilde M$ along 
$\gamma_{v^L_u}$ is given by 
$$FR_{v^L_u}=\{\frac{1+kb_j}{\bar{\lambda}_j(v^L_u)}\,\vert\,j=1,\cdots,
r,\,\,k\in{\bf Z}\}.$$
Since $v^L$ is parallel with respect to the normal connection of 
$\widetilde M$, the value $\bar{\lambda}_j(v^L_u)$ is independent of the 
choice of $u\in\widetilde M$.  
Hence we denote this value by $c^v_j$.  
Assume that $v$ satisfies the following condition:

\vspace{0.2truecm}

($\ast_1$) $\frac{1+kb_j}{c_j^v}$ 
($j=1,\cdots,r,\,\,k\in{\bf Z}$) are mutually distinct.

\vspace{0.2truecm}

\noindent
Note that $\{v_x\,\vert\,v\,\,{\rm satisfies}\,\,(\ast_1)\}$ is dense 
in $T^{\perp}_xM$.  Since $\pi\circ\phi$ is a Riemannian submersion, 
we have $FR_{v_x}=FR_{v^L_u}$ and hence 
${\rm Ker}\,Q_{v_x}(\frac{1+kb_j}{c_j^v})\not=\{0\}$ 
($j=1,\cdots,r,\,\,k\in{\bf Z}$).  Define a distribution $\overline E_{jk}$ 
on $M$ by $(\overline E_{jk})_x:=
\displaystyle{{\rm Ker}\,Q_{v_x}\left(\frac{1+kb_j}{c_j^v}\right)}$ and 
a distribution $E_{jk}$ on $\widetilde M$ by 
$(E_{jk})_u:=\displaystyle{{\rm Ker}\left(\widetilde A^{\bf c}_{v^L_u}-
\frac{c_j^v}{1+kb_j}{\rm id}\right)}$.  It is easy to show 
that $(\pi\circ\phi)^{\bf c}_{\ast u}((E_{jk})_u)
=(\overline E_{jk})_{(\pi\circ\phi)(u)}$ for any $u\in\widetilde M$.  
Fix $x\in M$ and $u\in(\pi\circ\phi)^{-1}(x)$.  Also, we can show that 
$\overline E_{jk}$'s ($k\in{\bf Z}$) coincide with one another.  
Hence, for each $j\in\{1,\cdots,r\}$, we have 
$$
Q_{v_x}(\frac{1+kb_j}{c_j^v})\vert_{(\overline E_{j0})_x}
=\left(\cosh(\frac{1+kb_j}{c_j^v}\sqrt{R(v_x)})-
\frac{\sinh(\frac{1+kb_j}{c_j^v}\sqrt{R(v_x)})}
{\sqrt{R(v_x)}}\circ A_{v_x}^{\bf c}\right)\vert_{(\bar E_{j0})_x}=0
$$
($k\in{\bf Z}$).  Set $D^x_{\beta}:={\rm Ker}(R(v_x)-\beta\,{\rm id})$ 
($\beta\in{\rm Spec}\,R(v_x)$).  
Fix $X(\not=0)\in(\bar E_{j0})_x$.  
Let $S_x$ be the minimal subset of ${\rm Spec}\,R(v_x)$ 
satisfying 
$\{X,A_{v_x}^{\bf c}X\}\subset\displaystyle{\mathop{\oplus}_{\beta\in S_x}
D_{\beta}^x}$.  Then, according to Lemma 3.2, we have 
$\displaystyle{\frac{b_j\sqrt{-\beta}}{c_j^v}\in\pi\sqrt{-1}{\bf Z}}$ for any 
$\beta\in S_x$.  Since $b_j$ and $c_j^v$ are is independent 
of the choice of $x\in M$, this fact together with the arbitrarinesses of 
$X$ and $j$ that ${\rm Spec}R(v_x)$ is independent of the choice of $x\in M$.  

\noindent
{\bf (Step II)} 
Take a parallel normal vector field $v$ of $M$.  
Let $\eta_{sv}(:M\to G/K)$ be the end-point map for $sv$ and 
$M_{sv}:=\eta_{sv}(M)$, where $s$ is sufficiently close to zero.  
Define a function $F_{sv}$ on $M$ by $\eta_{sv}^{\ast}\omega_{sv}
=F_{sv}\omega$, where $\omega$ (resp. $\omega_{sv}$) is the volume element of 
$M$ (resp. $M_{sv}$).  
Set $F_{v_x}(s):=F_{sv}(x)$ ($x\in M$), which coincides with the function 
$F_{v_x}$ stated in Introduction.  It is shown that 
$F_{v_x}$ has the holomorphic extension $F^{\bf c}_{v_x}$.  
According to Corollary 2.6 of [HLO], $M$ is isoparametric if and only if 
the projection from $M$ to any (sufficiently close) parallel submanifold along 
the sections is volume preserving up to a constant factor.  
Hence we suffice to show that $F_{v_{x_1}}^{\bf c}=F_{v_{x_2}}^{\bf c}$ holds 
for any two points $x_1$ and $x_2$ of $M$ in order to show that $M$ is 
isoparametric.  
By the way, since the complex focal radii along the geodesic $\gamma_{v_x}$ 
occur as zero points of $F_{v_x}^{\bf c}$, it follows from the complex 
equifocality of $M$ that $(F_{v_{x_1}}^{\bf c})^{-1}(0)
=(F_{v_{x_2}}^{\bf c})^{-1}(0)$ holds for any two points $x_1$ and $x_2$ of 
$M$.  Take a continuous orthonormal tangent frame field $(e_1,\cdots,e_n)$ of 
$M$ defined on a connected open set $U$ such that 
$R(v)(e_i)=\beta_ie_i$ ($i=1,\cdots,n$).  By the fact shown in Step I, 
$\beta_i$ ($i=1,\cdots,n$) are constant on $U$.  
Let $A_ve_i=\sum\limits_{j=1}^na_{ij}e_j$ ($i=1,\cdots,n$), 
where $a_{ij}$ ($i,j=1,\cdots,n$) are continuous functions on $U$.  
The strongly $M$-Jacobi field $J_{i,x}$ 
along $\gamma_{v_x}$ ($x\in U$) with $J_{i,x}(0)=e_{ix}$ is described as 
$$J_{i,x}(s)=\sum_{j=1}^n\left(\cosh(s\sqrt{-\beta_i})\delta_{ij}-
\frac{a_{ij}(x)\sinh(s\sqrt{-\beta_j})}
{\sqrt{-\beta_j}}\right)P_{\gamma_{v_x}\vert_{[0,s]}}e_{jx},$$
where $\delta_{ij}$ is the Kronecker's symbol and 
$\frac{\sinh(s\sqrt{-\beta_j})}{\sqrt{-\beta_j}}$ implies $s$ when 
$\beta_j=0$.  From this description, we have 
$$F_{v_x}(s)={{\rm det}}\left(\cosh(s\sqrt{-\beta_i})\delta_{ij}
-\frac{a_{ij}(x)\sinh(s\sqrt{-\beta_j})}{\sqrt{-\beta_j}}\right),$$
where $\left(\cosh(s\sqrt{-\beta_i})\delta_{ij}-\frac{a_{ij}(x)
\sinh(s\sqrt{-\beta_j})}{\sqrt{-\beta_j}}\right)$ 
is the matrix of $(n,n)$-type whose $(i,j)$-component is 
$\cosh(s\sqrt{-\beta_i})\delta_{ij}-\frac{a_{ij}(x)\sinh(s\sqrt{-\beta_j})}
{\sqrt{-\beta_j}}$.  
Hence we have 
$$F_{v_x}^{\bf c}(z)={{\rm det}}
\left(\cos(\sqrt{-1}z\sqrt{-\beta_i})\delta_{ij}-\frac{a_{ij}(x)
\sin(\sqrt{-1}z\sqrt{-\beta_j})}{\sqrt{-1}\sqrt{-\beta_j}}\right). 
\leqno{(3.5)}$$
Define a subset ${\cal W}$ of $T^{\perp}_xM$ by 
${\cal W}:=\{w\in T^{\perp}_xM\,\vert\,\sqrt{-b_1^w}:\cdots:\sqrt{-b^w_n}\,\,{\rm is}\,\,{\rm an}\,\,
{\rm integer}\,\,{\rm ratio}\}$, where $\{b^w_1,\cdots,b^w_n\}$ 
($b_1^w\leq\cdots\leq b^w_n$) is all eigenvalues of $R(w)$.  
Let $x=gK$ and $\mathfrak a:=g_{\ast}^{-1}T^{\perp}_xM$.  Since 
${\rm codim}\,M={\rm rank}(G/K)$ by the assumption, $\mathfrak a$ is a maximal 
abelian subspace of $\mathfrak p:=T_{eK}(G/K)\,(\subset\mathfrak g)$.   Let 
$\triangle$ be the root system of $G/K$ with respect to $\mathfrak a$.  Since 
the set of all the eigenvalues of $R(w)$ is equal to 
$\{-\alpha(w)^2\,\vert\,\alpha\in\triangle\}$, ${\cal W}$ is dense in 
$T_x^{\perp}M$.  Assume that $v_x$ belongs to ${\cal W}$.  
Then, since $\sqrt{-\beta_1(x)}:\cdots:\sqrt{-\beta_n(x)}$ is an integer 
ratio, they are expressed as $\sqrt{-\beta_i(x)}=m_ib$ ($i=1,\cdots,n$) 
in terms of some real constant $b$ and integers $m_1,\cdots,m_n$.  
Hence the function $F^{\bf c}_{v_x}$ ($x\in U$) is described as 
$F^{\bf c}_{v_x}(z)=e^{-bz\sum_{i=1}^n\vert m_i\vert}G_x(e^{2bz})$ in terms of 
some polynomial $G_x$ of degree $\vert m_1\vert+\cdots+\vert m_n\vert$.  
Take arbitrary two points $x_1$ and $x_2$ of $U$.  Let $c:[0,1]\to U$ be 
a continuous curve with $c(0)=x_1$ and $c(1)=x_2$.  Since $M$ is complex 
equifocal, $(F^{\bf c}_{v_{c(t)}})^{-1}(0)$ is independent of the choice of 
$t\in [0,1]$.  Hence so is also $G_{c(t)}^{-1}(0)$.  This implies that 
$G_{x_1}=aG_{x_2}$ for some a non-zero complex constant $a$.  
Hence we have $F^{\bf c}_{x_1}=aF^{\bf c}_{x_2}$.  
Furthermore, since $F^{\bf c}_{v_{x_1}}(0)=F^{\bf c}_{v_{x_2}}(0)=1$, we have 
$F^{\bf c}_{v_{x_1}}=F^{\bf c}_{v_{x_2}}$.  From the arbitarinesses of 
$x_1,\,x_2$ and $U$, we see that $F^{\bf c}_{v_x}$ is independent of the 
choice of $x\in M$.  Furthermore, since ${\cal W}$ is dense in $T^{\perp}_xM$, 
it follows that $F^{\bf c}_{v_x}$ is independent of the choice of $x\in M$ 
in the case where $v$ does not belong to ${\cal W}$.  
Thus $M$ is an isoparametric submanifold with flat section.  
\hspace{1.5truecm}q.e.d.

\vspace{0.5truecm}

Let an abelian subspace $\mathfrak b$ of $\mathfrak p:=T_{eK}(G/K)$, 
$\mathfrak a$ be a maximal abelian subspace of $\mathfrak p$ containing 
$\mathfrak b$ and $\triangle$ be the root system with respect to 
$\mathfrak a$.  Set $\triangle_{\mathfrak b}:=\{\alpha\vert_{\mathfrak b}\,
\vert\,\alpha\in\triangle\,\,{\rm s.t.}\,\,\alpha\vert_{\mathfrak b}\not=0\}$, 
where $\alpha\vert_{\mathfrak b}$ is the restriction of $\alpha$ to 
$\mathfrak b$.  This set $\triangle_{\mathfrak b}$ is independent of the 
choice of $\mathfrak a$.  
By using Lemma 3.2, we prove the following lemma.  

\vspace{0.5truecm}

\noindent
{\bf Lemma 3.3.} {\sl Let $M$ be a proper complex equifocal submanifold in 
a symmetric space $G/K$ of non-compact type.  
If, for any $\beta_1,\,\beta_2\in\triangle_{\mathfrak b}$, $\beta_1$ and 
$\beta_2$ are linearly independent (over ${\Bbb R}$) or 
$\beta_1=\pm\beta_2$, then $M$ is curvature-adapted.}

\vspace{0.5truecm}

\noindent
{\it Proof.} 
Let $\{\lambda_i\,\vert\,i\in I\}$ be the set of all complex principal 
curvatures of $\widetilde M$ at $\hat 0$ and set 
${\it l}_i:=\lambda_i^{-1}(1)$ ($i\in I$), wehre $\widetilde M$ is 
as above.  Since the group generated by the complex reflections $R_i$'s 
($i\in I$) of order $2$ with respect to ${\it l}_i$ is discrete, 
${\cal L}:=\{{\it l}_i\,\vert\,i\in I\}$ is equal to the sum 
$\displaystyle{\mathop{\cup}_{j=1}^r{\cal L}_j}$ of subfamilies 
${\cal L}_j:=\{{\it l}_{i^j_k}\,\vert\,k\in{\bf Z}\}$ ($j=1,\cdots,r$) of 
parallel complex hyperplanes equidistant to one another.  
For each $j\in\{1,\cdots,r\}$, we have $\lambda_{i^j_k}^{-1}(1)
=\lambda_{i^j_0}^{-1}(1+kb_j)$ for some $b_j\in{\bf C}$.  For simplicity, 
we denote $\lambda_{i^j_0}$ by $\overline{\lambda}_j$ $(j=1,\cdots,r$).  
Let $\mathfrak a$ be a maximal abelian subspace of 
$\mathfrak p:=T_{eK}(G/K)$ containing $\mathfrak b:=T_{eK}^{\perp}M$, 
$\triangle$ be the root system with respect to $\mathfrak a$ and 
$\mathfrak p=\mathfrak a+\sum\limits_{\alpha\in\triangle_+}
\mathfrak p_{\alpha}$ be the root space decomposition with respect to 
$\mathfrak a$.  
Set $U:=\{v\in\mathfrak b\,\vert\,\beta(v){\rm 's}\,\,
(\beta\in\triangle_{\mathfrak b})\,\,{\rm are}\,\,{\rm linearly}\,\,
{\rm independent}\,\,{\rm over}\,\,{\bf Q}\}$, which is dense in 
$\mathfrak b$ by the assumption for $\triangle_{\mathfrak b}$.  
Fix $v\in U$.  
The set $FR_v$ of all complex focal radii of $M$ along $\gamma_v$ is given by 
$$FR_v=\{z\in{\bf C}\,\vert\,{\rm Ker}\,Q_v(z)\not=\{0\}\},$$
where $Q_v(z):=D^{co}_{zv}-D^{si}_{zv}\circ A^{\bf c}_{zv}$.  
The constant path $\hat v$ at $v$ is the horizontal lift of $v$ to $\hat 0$.  
On the other hand, the 
set $FR_{\hat v}$ of all complex focal radii of $\widetilde M$ along 
$\gamma_{\hat v}$ is given by 
$$FR_{\hat v}=\{\frac{1+kb_j}{\bar{\lambda}_j(v)}\,\vert\,j=1,\cdots,r,\,\,
k\in{\bf Z}\}.$$
Since $FR_v=FR_{\hat v}$, we have 
${\rm Ker}\,Q_v(\frac{1+kb_j}{\bar{\lambda}_j(v)})\not=\{0\}$ 
($j=1,\cdots,r,\,\,k\in{\bf Z}$).  Set $\overline E_{jk}^v:=
\displaystyle{{\rm Ker}\,Q_v\left(\frac{1+kb_j}{\bar{\lambda}_j(v)}\right)}$ 
and 
$E_{jk}^v:=\displaystyle{{\rm Ker}\left(\widetilde A^{\bf c}_{\hat v}-
\frac{\bar{\lambda}_j(v)}{1+kb_j}{\rm id}\right)}$.  It is easy to show that 
$(\pi\circ\phi)^{\bf c}_{\ast}(E_{jk}^v)={\overline E}_{jk}^v$.  Also, we 
can show that ${\overline E}_{jk}^v$'s ($k\in{\bf Z}$) coincide with one 
another.  Hence, for each $j\in\{1,\cdots,r\}$, we have 
$$
Q_v(\frac{1+kb_j}{\bar{\lambda}_j(v)})\vert_{{\overline E}_{j0}^v}
=\left(\cosh(\frac{1+kb_j}{\bar{\lambda}_j(v)}{\rm ad}(v))-
\frac{\sinh(\frac{1+kb_j}{\bar{\lambda}_j(v)}{\rm ad}(v))}
{{\rm ad}(v)}\circ A_v^{\bf c}\right)\vert_{{\bar E}_{j0}^v}=0
$$
($k\in{\bf Z}$).  Fix $X(\not=0)\in\bar E_{j0}$.  
Set $\mathfrak p_{\beta}:=\displaystyle{
\mathop{\oplus}_{\alpha\in\triangle_+\,\,{\rm s.t.}\,\,
\alpha\vert_{\mathfrak b}=\pm\beta}\mathfrak p_{\alpha}}$ 
($\beta\in\triangle_{\mathfrak b}$).  Let $S$ be the minimal 
subset of $\triangle_{\mathfrak b}$ satisfying 
$\{X,A_v^{\bf c}X\}\subset\displaystyle{\mathop{\oplus}_{\beta\in S}
\mathfrak p_{\beta}^{\bf c}}$.  Then, according to Lemma 3.2, we have 
$\displaystyle{\frac{b_j\beta(v)}{\overline{\lambda}_j(v)}\in\pi\sqrt{-1}
{\bf Z}}$ for any $\beta\in S$.  
Hence, since $\beta(v)$'s ($\beta\in\triangle_{\mathfrak b}$) are linearly 
independent over ${\Bbb Q}$, we see that $S$ is a one-point set.  
That is, we have 
$X\in\mathfrak p_{\beta_0}^{\bf c}$ and 
$A_v^{\bf c}X\in\mathfrak p_{\beta_0}^{\bf c}$ for some 
$\beta_0\in\triangle_{\mathfrak b}$.  
Hence we have $A_v^{\bf c}X=\frac{\beta_0(v)}
{\tanh({\beta_0(v)}/\bar{\lambda}_j(v))}X$.  
This together with the arbitrariness of 
$X(\in{\overline E}^v_{j0})$ implies 
$R(v)^{\bf c}({\overline E}^v_{j0})\subset{\overline E}^v_{j0}$ and 
$[A_v^{\bf c},R(v)^{\bf c}]\vert_{{\bar E}^v_{j0}}=0$.  
On the other hand, since $M$ is proper complex equifocal, we have 
$$\left(\mathop{\oplus}_{j=1}^r{\bar E}_{j0}^v\right)\oplus
\left({\rm Ker}\,A_v^{\bf c}\cap{\rm Ker}\,R(v)^{\bf c}\right)
=(T_{eK}M)^{\bf c}.$$
Therefore we have 
$R(v)^{\bf c}((T_{eK}M)^{\bf c})\subset (T_{eK}M)^{\bf c}$ and 
$[A_v^{\bf c},R(v)^{\bf c}]=0$.  Hence we have 
$R(v)(T_{eK}M)\subset T_{eK}M$ and $[A_v,R(v)]=0$.  
Furthermore, since $v$ is an arbitrary element of $U$ and 
$U$ is dense in $\mathfrak b$, 
$R(w)(T_{eK}M)\subset T_{eK}M$ and $[A_w,R(w)]=0$ holds 
for any $w\in \mathfrak b$.  By the same 
discussion, we can show that the same fact holds for any point of $M$ other 
than $eK$.  Therefore $M$ is curvature-adapted.  \hspace{1.5truecm}q.e.d.

\vspace{0.5truecm}

By using Lemma 3.3., we prove Theorem C.  

\vspace{0.5truecm}

\noindent
{\it Proof of Theorem C.} From ${\rm codim}\,M={\rm rank}(G/K)$, we have 
$\triangle_{\mathfrak b}=\triangle$, where $\triangle$ and 
$\triangle_{\mathfrak b}$ are as above.  
Since $\triangle$ is reduced by the assumption, 
$\triangle_{\mathfrak b}(=\triangle)$ satisfies the condition in Lemma 3.3.  
Hence, it follows from Lemma 3.3 that $M$ is curvature-adapted.  
\hspace{0.2truecm}q.e.d.

\newpage


\centerline{{\bf References}}

\vspace{0.5truecm}

{\small

\noindent
[B1] J. Berndt, Real hypersurfaces with constant principal curvatures 
in complex hyperbolic 

space, J. Reine Angew. Math. {\bf 395} (1989) 132-141. 

\noindent
[B2] J. Berndt, Real hypersurfaces in quaternionic space forms, 
J. Reine Angew. Math. {\bf 419} 

(1991) 9-26. 



\noindent
[BT] J. Berndt and H. Tamaru, Cohomogeneity one actions on noncompact 
symmetric spaces 

with a totally geodesic singular orbit, 
Tohoku Math. J. {\bf 56} (2004) 163-177.

\noindent
{\small [BV] J. Berndt and L. Vanhecke, 
Curvature adapted submanifolds, 
Nihonkai Math. J. {\bf 3} (1992) 

177-185.



\noindent
[Ch] U. Christ, 
Homogeneity of equifocal submanifolds, J. Differential Geometry 
{\bf 62} (2002) 1-15.

\noindent
[GT] O. Goertsches and G. Thorbergsson, 
On the Geometry of the orbits of Hermann actions, 

Geom. Dedicata {\bf 129} (2007) 101-118.

\noindent
[HLO] E. Heintze, X. Liu and C. Olmos, 
Isoparametric submanifolds and a 
Chevalley type rest-

riction theorem, Integrable systems, geometry, and topology, 151-190, 
AMS/IP Stud. Adv. 

Math. 36, Amer. Math. Soc., Providence, RI, 2006.

\noindent
[HPTT] E. Heintze, R.S. Palais, C.L. Terng and G. Thorbergsson, 
Hyperpolar actions on symme-

tric spaces, Geometry, topology and physics for Raoul Bott (ed. S. T. Yau), Conf. Proc. 

Lecture Notes Geom. Topology {\bf 4}, 
Internat. Press, Cambridge, MA, 1995 pp214-245.

\noindent
[He] S. Helgason, 
Differential geometry, Lie groups and symmetric spaces, 
Academic Press, New 

York, 1978.

\noindent
[Koi1] N. Koike, 
Proper isoparametric semi-Riemannian submanifolds in a 
semi-Riemannian 

space form, Tsukuba J. Math. {\bf 13} (1989) 131-146.

\noindent
[Koi2] N. Koike, 
Submanifold geometries in a symmetric space of non-compact 
type and a pseudo-

Hilbert space, Kyushu J. Math. {\bf 58} (2004) 167-202.

\noindent
[Koi3] N. Koike, 
Complex equifocal submanifolds and infinite dimensional anti-
Kaehlerian isopara-

metric submanifolds, Tokyo J. Math. {\bf 28} (2005) 201-247.

\noindent
[Koi4] N. Koike, 
Actions of Hermann type and proper complex equifocal submanifolds, 
Osaka J. 

Math. {\bf 42} (2005) 599-611.

\noindent
[Koi5] N. Koike, 
A splitting theorem for proper complex equifocal submanifolds, Tohoku Math. J. 

{\bf 58} (2006) 393-417.

\noindent
[Koi6] N. Koike, Complex hyperpolar actions with a totally geodesic orbit, 
Osaka J. Math. {\bf 44} 

(2007) 491-503.

\noindent
[Koi7] N. Koike, 
The homogeneous slice theorem for the complete complexification of a 
proper 

complex equifocal submanifold, Tokyo J. Math. {\bf 33} (2010), 1-30.

\noindent
[Koi8] N. Koike, Examples of a complex hyperpolar action without 
singular orbit, Cubo A Math. 

J. {\bf 12} (2010), 127-144.  

\noindent
[Koi9] N. Koike, Hermann type actions on a pseudo-Riemannian symmetric 
space, Tsukuba J. 

Math. (to appear) (arXiv:math.DG/0807.1604v2).

\noindent
[Koi10] N. Koike, Homogeneity of irreducible proper complex equifocal 
submanifolds, arXiv:math.

DG/0807.1606v2.

\noindent
[Koi11] N. Koike, On curvature-adapted and proper complex equifocal 
submanifolds, Kyungpook 

Math. J. (to appear) (arxiv:math.DG/0809.4933v1).

\noindent
[Kol] A. Kollross, A Classification of hyperpolar and cohomogeneity one 
actions, Trans. Amer. 

Math. Soc. {\bf 354} (2001) 571-612.


\noindent
[PT] R.S. Palais and C.L. Terng, Critical point theory and submanifold 
geometry, Lecture Notes 

in Math. {\bf 1353}, Springer, Berlin, 1988.

\noindent
[T1] C.L. Terng, 
Isoparametric submanifolds and their Coxeter groups, 
J. Differential Geometry 

{\bf 21} (1985) 79-107.

\noindent
[T2] C.L. Terng, 
Proper Fredholm submanifolds of Hilbert space, 
J. Differential Geometry {\bf 29} 

(1989) 9-47.


\noindent
[TT] C.L. Terng and G. Thorbergsson, 
Submanifold geometry in symmetric spaces, 
J. Differential 

Geometry {\bf 42} (1995) 665-718.

\noindent
[W] B. Wu, Isoparametric submanifolds of hyperbolic spaces, 
Trans. Amer. Math. Soc. {\bf 331} 

(1992) 609-626.
}

\vspace{1truecm}

\rightline{Department of Mathematics, Faculty of Science, }
\rightline{Tokyo University of Science}
\rightline{1-3 Kagurazaka Shinjuku-ku,}
\rightline{Tokyo 162-8601, Japan}
\rightline{(e-mail: koike@ma.kagu.tus.ac.jp)}

\end{document}